\documentclass[12pt,a4paper]{article} 
\usepackage{amsfonts,amssymb} 

\usepackage{tikz}

  \def\?#1{}
\def\IZ{{\mathbb{Z}}}  
   \def\cB{{\cal B}} 
\def\cT{{\cal T}}   \def\cF{{\cal F}}

\def\ep{\varepsilon}   \def\bdelete#1{}

\def\beq#1#2{\begin{equation} \label{#1} #2 \end{equation}}
 
\def\beaq#1#2{\label{#1} \begin{eqnarray} #2 \end{eqnarray}}

\def\function#1{\left\{\!\!\!\begin{array}{ll} #1 \end{array} \right.}

\def\proof{\smallskip \noindent {\bf Proof. \ }}       
\def\blanksquare{\,\,\,$\sqcup\!\!\!\!\sqcap$}         
\def\qed{\hfill\blanksquare\linebreak\smallskip\par}   

\def\thname{Theorem}  \def\lmname{Lemma}    \def\prname{Proposition}
\def\dfname{Definition}  \def\crname{Corollary}  \def\rmname{Remark}
\def\exname{Example}

\newtheorem{theorem}{\thname}
\newtheorem{lemma}{\lmname}
\newtheorem{corollary}[lemma]{\crname}   
\newtheorem{example}{\exname}

\newtheorem{dftn}{\dfname}
\newenvironment{definition}{\begin{dftn}\rm}{\end{dftn}} 
\def\bdef#1{\begin{definition} #1 \end{definition}}
\newtheorem{rmrk}[lemma]{\rmname}

\makeatletter\def\fps@figure{htbp}\makeatother 

\begin{document}

\title{
        Topological and metric recurrence \\ for general Markov chains}
\author{Michael Blank\thanks{
        Institute for Information Transmission Problems RAS
        (Kharkevich Institute);}
        \thanks{National Research University ``Higher School of Economics'';
        e-mail: blank@iitp.ru}
       \\ \\{\small The paper is dedicated to the memory of Robert Adolfovich Minlos.}
       }
\date{June 1, 2018} 
\maketitle

\begin{abstract}%
Using ideas borrowed from topological dynamics and ergodic theory
we introduce topological and metric versions of the recurrence
property for general Markov chains. The main question of interest
here is how large is the set of recurrent points. We show that 
under some mild technical assumptions the set of non recurrent
points is of zero reference measure. Necessary and sufficient 
conditions for a reference measure $m$ (which needs not to be 
dynamically invariant) to satisfy this property are obtained. 
These results are new even in the purely deterministic setting.
\end{abstract}%

\section{Introduction}\label{s:intro}

The aim of this work is to analyze a circle of questions related
to the notion of recurrence in general Markov chains. Being well
known in two very different subfields of random systems: lattice
random walks and ergodic theory of continuous selfmaps, the
recurrence property is next to being neglected in general theory
of Markov chains (perhaps except for a few notable exceptions
which we will discuss in detail).

There are two known approaches to the definition of a recurrent
point: topological and metric (defined in very different
contexts). In order to introduce their analogies for general
Markov chains let us recall basic definitions.

\bdef{By an inhomogeneous Markov chain one means a random process
$\xi_t: (\Omega,\cF,P) \to (X,\cB, m)$ acting on a Borel $(X,\cB)$
space with a finite reference measure $m$ (which needs not to
coincide with the distribution of the process $\xi_t$). This
process is completely defined by a family of {\em transition
probabilities}
$$Q_s^t(x,A):=P(\xi_{s+t}\in A|\xi_s=x),~A\in\cB.$$
If the chain is homogeneous, i.e. the transition probabilities do
not depend on $s$ we drop the lower index and write
$Q^t(x,A)\equiv Q_s^t(x,A)$. }

\bdef{By the $t$-{\it preimage} with $t\ge0$ of a set $B\in\cB$
under the action of the homogeneous Markov chain $\xi_t$ we call
the set of points
$$ Q^{-t}(B) := \{x\in X:~~ Q^t(x,B)>0\} .$$
}
In other words this is the set of initial points of trajectories which 
reach the set $B$ at time $t$ with positive probability.  

Now we are ready to return to the notion of recurrence. Observe
that since the phase space is equipped with the Borel
$\sigma$-algebra, it is equipped with the corresponding topology
as well. We start from the notion of the topological recurrence
well known in the field of topological dynamics (see e.g.
\cite{Ni}).

\bdef{ A point $x\in X$ is called {\em topologically recurrent} if
for any open neighborhood $U\ni x$ for each $s$ there exists an
(arbitrary large) $t=t(x,U,s)$ such that $Q_s^t(x,U)>0$ (i.e. a
trajectory eventually returns to $U$ with positive probability). }

\bdef{ A point $x\in X$ is called {\em metrically recurrent} if
for any set $V\ni x$ of positive $m$-measure and any $s$ there
exists an (arbitrary large) $t=t(x,V,s)$ such that $Q_s^t(x,V)>0$
(i.e. a trajectory eventually returns to $V$ with positive
probability). }

The last definition is our modification of the metric version of
the recurrence property proposed by T.E. Harris in \cite{Ha} in
order to get reasonably general assumptions guaranteeing the
existence of an invariant measure. In fact, Harris used this
property only in the case when the reference measure $m$ is
invariant with respect to the process. Another weak point of
Harris' approach is that whence a point $x$ is metrically
recurrent, the corresponding trajectory (realization of the
process) needs to visit any set of positive measure with positive
probability, which looks way too excessive.

The 3d approach to the recurrence notion is related to ergodic 
theory of deterministic dynamical systems, where it is well 
known and studied, but again only in the case when the 
measure $m$ is dynamically invariant. The probabilistic version 
may be formulated as follows. 

\bdef{ A point $x\in X$ is called {\em Poincare recurrent} with
respect to a $\cB$-measurable set $A\ni x$ if for each $s$ there
exists an (arbitrary large) $t=t(x,A,s)$ such that $Q_s^t(x,A)>0$
(i.e. a trajectory eventually returns to $A$ with positive
probability). }

Comparing these definitions with their deterministic counterparts
(see e.g. \cite{Ni,Si}) or to the notions of recurrence and
transience well studied for the case of countable Markov chains
(see e.g. \cite{Fe,MT,Or,Sp}) one is tempted to make the
conditions stronger assuming that the corresponding events take
place with probability one (instead of just being positive).
Unfortunately, as we will show in section~\ref{s:gen} the results
we are looking for do not hold under these stronger assumptions.
To some extent this may explain why these properties were not
studied earlier.

The main question of interest for us is to find how large is the
set of recurrent points?

The celebrated Poincare recurrence theorem (see e.g. \cite{Si})
states that for a measurable discrete time dynamical system
$(T,X)$ and any measurable set $A$ the set of non Poincare
recurrent points is of zero measure with respect to any
$T$-invariant measure. There is a number of obstacles which one
needs to overcome to generalize this result for a general Markov
chain: continuous time, inhomogeneity, and more important a
principal absence of invariant measures (stationary distributions)
for inhomogeneous Markov chains.\footnote{Note that non
translation-invariant infinite volume  Gibbs measures might be 
still present in this setting.} 
Here we formulate our main result
in this direction only for the discrete time case ($t\in\IZ$)
leaving the more complicated continuous time version for
Section~\ref{s:gen}.


\begin{theorem}\label{t:main}
Let $m$ be a finite measure on $(X,\cB)$ and let $Q_s^t$ does not
depend on $s$. Then the property that for each set $A\in\cB$ the
set of Poincare recurrent points in $A$ is of full $m$-measure
(i.e. its complement in $A$ is of zero $m$-measure)
is equivalent to: %
\beaq{e:main}{\sum_{n\ge1} m(Q^{-n}(A)\cap A) = \infty ~~~\forall
A\in\cB:~m(A)>0. }
\end{theorem}

An important observation here is that if $Q_s^t$ does depend on
$s$, then the direct statement in Theorem~\ref{t:main} (that (\ref{e:main}) 
implies the abundance of Poincare recurrent points) remains
correct, but the inverse one may fail (see Section~\ref{s:dis}).

The situation with other types of recurrence is much more subtle,
in particular a very ``good'' Markov chain, for which all points 
are topologically recurrent with respect to the ``standard'' 
topology, may change drastically when one chooses a
different topology instead. Under this new topology all points
might become non-recurrent, e.g. irrational circle rotation 
with the discrete topology.

The following result gives conditions for the topological
recurrence of ``typical'' points. Here and in the sequel 
we assume that the measurable space $(X,\cB)$ is equipped 
with a topology $\cT$, which is compatible with the 
$\sigma$-algebra $\cB$.

\begin{theorem} \label{t:recur-t} Let $(X,\cT)$ be a compact 
metric space and the finite reference measure $m$ satisfy the 
property~(\ref{e:main}). Then $m$-almost every point $x\in X$ is
topologically and metrically recurrent.
\end{theorem}

It is natural to ask what can we say about recurrent points 
if the property~(\ref{e:main}) does not hold. As we already noted 
if there is a dynamical invariant measure $\mu$ then one can 
easily construct from it a reference measure $m$ satisfying 
(\ref{e:main}). For example, any finite measure absolutely 
continuous with respect to $\mu$. Therefore to answer to this 
question one needs to consider systems without invariant measures. 
This will be done in section~\ref{s:gen}, where among other results 
we will demonstrate that a measurable Markov chain with a compact 
phase space may have no recurrent points. 

\bigskip

It is worth noting that the results formulated in
Theorems~\ref{e:main} and \ref{t:recur-t} are new even in the
purely deterministic setting.

The paper is organized as follows. In the next section we review a
few definitions related to Markov chains which are necessary for
the further analysis. The sections \ref{s:proof-T1} and
\ref{s:proof-T2} are dedicated to the proof of
Theorems~\ref{t:main} and \ref{t:recur-t} respectively. Then in
section~\ref{s:dis} we study connections between the condition
(\ref{e:main}) and other ``conservativity'' type properties. One
of interesting moments here is a striking difference between
forward and backward in time statistics. Finally in
section~\ref{s:gen} we discuss some generalizatins. Namely we 
show that a stronger version of the recurrence property plays 
a very different role than the one defined above; discuss 
recurrence in the absence of invariant measures; give continuous 
time versions of our main results; and finally analyze their connection 
to yet to be proven stochastic version of the so called multiple 
recurrence problem.

\section{Preliminaries}\label{s:prelim}

The transition probabilities $Q_s^t(\cdot,\cdot)$ satisfy the
following standard properties:
\begin{itemize}
\item For fixed $s,t,x$ the function $Q_s^t(x,\cdot)$ is a
probability measure on the $\sigma$-algebra $\cB$.

\item For fixed $s,t,A$ the function $Q_s^t(\cdot,A)$ is
$\cB$-measurable.

\item For $t=0$ ~~$Q_s^t(x,A)=\delta_x(A)$.

\item For each $s, 0\le t \le t'$ and $A\in\cB$ we have
$$  Q_s^{t'}(x,A) = \int_X Q_s^t(x,dy) Q_t^{t'-t}(y,A) .$$

\end{itemize}

The process $\xi_t$ induces the action on measures:
$$Q_s^t\mu(A) := \int Q_s^t(x,A)~d\mu(x)$$
and the action on functions:
$$Q_s^t\phi(x) := \int \phi(y)Q_s^t(x,dy).$$
In particular, the well known Feller property in terms of the
action on functions means that $Q_s^t:C^0\to C^0~~\forall
s,t\ge0$.

\bdef{ A Borel measure $\mu$ is said to be {\em invariant} or {\em
stationary} for the Markov chain $\xi_t$ if it is a solution to
the equation
$$ Q_s^t\mu = \mu~~\forall s,t .$$
}

In the discrete time setting $t\in\IZ$ and we write $Q(\cdot,\cdot)\equiv Q^1(\cdot,\cdot)$. 

\bdef{We say that a discrete time homogeneous Markov chain 
on a measurable space $(X,\cB)$ defined by the transition probabilities 
$Q(\cdot,\cdot)$ is {\em deterministic} if there exists a measurable 
self-map $T:X\to X$ such that
$$ Q(x,A) := \function{1 &\mbox{if } Tx \in A \\
                       0 &\mbox{otherwise}.} $$
}

The deterministic nature of a Markov chain may be characterized
also in terms of the action on measures, namely that the image of
any Dirac measure is again a Dirac measure, i.e.
$$ \forall x\in X~~\exists y\in X:~~ Q\delta_x:=\delta_y .$$

By the properties of the transition probabilities, if $B\in\cB$,
then the $t$-preimage $Q^{-t}(B)$ is measurable as well $\forall
t\ge0$. In the deterministic case the map $T$ acts on measures as
follows:
$$ T\mu(A) := \mu(T^{-1}A) .$$
Therefore one might expect that a similar property holds true for
a general Markov chain, i.e.
$$ Q^t m(B) = m(Q^{-t}(B)) .$$

The following trivial example demonstrates that this is absolutely
not the case and it emphasizes that some care is necessary when
one is trying to apply arguments well known in the deterministic
dynamical systems theory for a general Markov chain.

\begin{example}\label{ex0} Let $\xi_n$ be a discrete time 2-state
Markov chain with the transition probability matrix $(q_{ij})$
with $0<q_{11}<q_{22}<1/2$, and let $A$ stand for the the 1st state.
\end{example}

Then $Q^{-1}(A)=X$, and for a given nontrivial distribution $m$ on
$X$ we have
$$ m(Q^{-1}(A))=m(X) ~~\ne~~ Qm(A)=q_{11}m(A) + q_{21}m(X\setminus A)
.$$

\bdef{We say that a Markov chain $\xi_t$ satisfies the {\em
Poincare recurrence property} with respect to a measure $\mu$ and
use the notation PRP($\mu)$ for this, if the set of Poincare
recurrent points in each set $A\in\cB$ is of full $m$-measure.}

In these terms Theorem~\ref{t:main} claims that the Poincare
recurrence property is equivalent to the property (\ref{e:main}).

\section{Proof of Theorem~\ref{t:main}}\label{s:proof-T1}

First we prove that the property (\ref{e:main}) implies that the
set of Poincare recurrent points is of full $m$-measure. For a
given set $A\in\cB$ with $m(A)>0$ denote by $\hat{A}$ the set of
non Poincare recurrent points in it, i.e.
$$ \hat{A}:= \{x\in A:~~ Q^t(x,A)=0~~\forall t>0\} .$$

Assume from the contrary that $m(\hat{A})>0$.

In the deterministic setting when the Markov chain is defined 
by a measurable self-map $T:X\to X$ one may use the assumption 
$$ \sum_{n\ge1} m(T^{-n}A) = \infty ,$$
which is even weaker than (\ref{e:main}), and argue that 
for any pair of different moments of time $t \ne s$ the 
corresponding preimages of the set $\hat{A}$ are disjoint. 
Indeed, for 
$$ x\in Q^{-t}(\hat{A}) \equiv T^{-t}(\hat{A}) $$ 
the forward ``orbit'' of $x$ visits the set $A\supseteq \hat{A}$ last
time at the moment $t$, i.e. $Q^s(x,A)=0$ if $t>s$.

Therefore
$$ m(X) \ge m(\cup_{n\ge1}Q^{-n}(\hat{A})) =
\sum_{n\ge1}m(Q^{-n}(\hat{A})) .$$

It remains to notice that the measure $m$ is finite, but the right
hand side is equal to infinity by the assumption (\ref{e:main}).
We came to the contradiction. 

Unfortunately for a general Markov chain the preimages of the set 
$\hat{A}$ need not to be disjoint. To demonstrate this consider 
an example.

\begin{example}\label{exM} Let $\xi_n$ be a discrete time 3-state 
Markov chain with the transition probability matrix
$$\left( \begin{array}{ccc}
                0 & 0 & 1 \\
                1/2 & 1/2 & 0 \\
                0 & 0 & 1 \end{array} \right) $$
and let the reference measure $m$ be uniform on the phase space.
\end{example}

The set $A:=\{1\}$, consisting of the 1st state of this Markov 
chain, is obviously non-recurrent, since this state returns 
back with probability zero. On the other hand, all its $t$-preimages 
for $t\ge1$ coincide with the 2nd state, hence they are not disjoint. 
Additionally, 
$$ \sum_{n\ge1}m(Q^{-n}(A)) = \sum_{n\ge1}m(Q^{-n}(\{2\})) = \infty ,$$
since the $m$-measure of the 2nd state is equal to $1/3$. 

\bigskip

Therefore to study a general Markov chain we need to use a 
slightly more complex assumption (\ref{e:main}) to deal with the 
analysis of the set of non Poincare recurrent points $\hat A$., 
for which 
$$ m(Q^{-n}(\hat{A}) \cap \hat A) \le m(Q^{-n}(\hat{A}) \cap A) = 0 .$$
Thus the property (\ref{e:main}) fails and we came to the 
contradiction, which implies that the set of non Poincare recurrent 
points $\hat{A}$ should have zero $m$-measure.

To prove the claim in the inverse direction we need to demonstrate
that the Poincare recurrence property implies the divergence of
the sum in (\ref{e:main}). Denote the partial sum
$$ S(N) := \sum_{n=1}^N m(Q^{-n}(A) \cap A) $$
and assume from the contrary that
$$ S:= \limsup_{N\to\infty} S(N) <\infty .$$
Then for each $\ep>0$ there exists $N_\ep<\infty$ such that %
\beaq{e:small-ass}{ 0\le S - S(N_\ep) < \ep .} %
On the other hand, by PRP($m$) for each $N<\infty$
$$ m(x\in A:~~Q^n(x,A)>0,~n\ge N) > m(A)/2 ,$$
which contradicts to (\ref{e:small-ass}) since all addends in
$S(N)$ are nonnegative.

Theorem~\ref{t:main} is proven. \qed


\section{Proof of Theorem~\ref{t:recur-t}}\label{s:proof-T2}

In this Section we again consider only homogeneous Markov chains.

In fact we will prove Theorem~\ref{t:recur-t} in a bit more
general setting, namely instead of compactness of the set $X$ we
assume only that the topological space $(X,\cT)$ has a countable
base $\{\beta_i\}_{i\in\IZ_+}$ and that $\mu$ is a finite
$\sigma$-additive measure on the measurable space $(X,\cB)$, where
the topology $\cT$ is compatible with the $\sigma$-algebra $\cB$.
This idea was proposed originally by my student A.~Zhevnerchuk in
the analysis of the purely deterministic dynamics with an
invariant measure $m$.\footnote{See also a close approach in 
http://planetmath.org/proofofpoincarerecurrencetheorem2 .} 
Here we use it for a general homogeneous
Markov dynamics and a reference measure $m$ which needs not to be
dynamically invariant.

\begin{lemma}\label{l:base}
Consider a family of binary valued measurable functionals
$\phi_B:X\to\{0,1\}$ indexed by measurable sets $B\in\cB$, such
that
$$\{x\in X:~~ \phi_B(x)=1\}\subseteq B ~~\forall B\in\cB .$$
Then $\forall B\in\cT$
$$ \mu(x\in X: \phi_{B}(x)=0) \le \sum_i \mu(x\in X: \phi_{\beta_i}(x)=0).$$
\end{lemma}

\proof For each open set $\beta$ from the base of topology 
introduce the set
$$ \alpha_i:= \{x\in X:~~ \phi_\beta(x)=0\} .$$
Then due to the definition of the topological base and the
$\sigma$-additivity of the measure $m$ we immediately have
$$ \mu(x\in X: \phi_{B}(x)=0) \le \sum_i \mu(\alpha_i).$$
\qed

Using this result we will prove Theorem~\ref{t:recur-t} as
follows. Let $x\in X$ and let $B$ be its any open neighborhood.
Since $B$ is a union of some elements from the countable base
$\{\beta_i\}$ there exists an element $\beta_{i(x)}$ such that
$x\in\beta_{i(x)}\subseteq B$.

Denote the functionals $\phi_A$ to be equal to $1$ if $Q^n(x,A)>0$
for some $n\in\IZ_+$ and to $0$ otherwise. Then the sets
$\alpha_i$ defined in the proof of Lemma~\ref{l:base} are the
subsets of the base sets $\beta_i$ whose points never return to
$\beta_i$ under dynamics. Thus $m(\alpha_i)=0~~\forall i$ by the
assumption on the PRP($m$). By Lemma~\ref{l:base} we get
$$ \mu(x\in X: \phi_{B}(x)=0) 
   \le \sum_i \mu(x\in X: \phi_{\beta_i}(x)=0) =0.$$
Therefore the measure of topologically non-recurrent points in $X$
is zero.

To prove that the measure of metrically non-recurrent points is
also zero one applies the argument above to a set of positive
$m$-measure instead of the open neighborhood.

Theorem~\ref{t:recur-t} is proven. \qed

\section{Discussion} \label{s:dis}

It is easy to see that whence the reference measure $m$ is
dynamically invariant or is absolutely continuous with respect to
an invariant measure, then the assumption~(\ref{e:main}) holds
true automatically. The aim of this section is to study the
connections of a somewhat unusual property~(\ref{e:main}) with
other assumptions of similar type, in particular with the
assumptions related to the conservativity property in dynamical
systems theory.

\subsection{Time inversion} \label{ss:time}

With respect to a given topology $\cT$ the definition of the 
$t$-preimage (with $t>0$) of a set $B$ may be reformulated as follows %
\beq{e:preim}{Q^{-t}(B) := \{y\in X:~~Q^t(y,U(B))>0~~
                             \forall U(B)\in\cT, ~B\subseteq U(B)\} .}
Similarly one defines a $t$-image having in mind the action 
on sets in the ``positive'' time direction.

\bdef{By the $t$-image with $t\ge0$ of a measurable set $A\in\cB$
we mean %
\beq{e:imag}{Q^{t}(A) := \cup_{x\in A}\{y\in X:~~Q^t(x,U(y))>0~~
                          \forall U(y)\in\cT, ~y\in U(y)\} .}
}

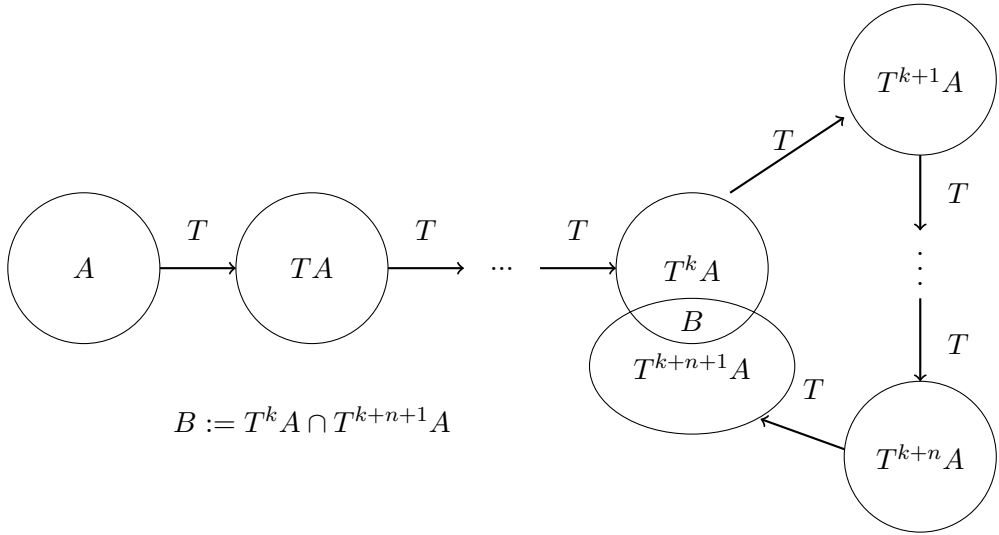
\begin{figure}
\begin{center}
\begin{tikzpicture} \small 
       \draw(0,0) circle (1cm); \node at (0,0) {$A$};  
       \draw  [thick, ->] (1,0) to (2,0); \node at (1.5,0.5) {$T$};
       \draw(3,0) circle (1cm); \node at (3,0) {$TA$};
       \draw  [thick, ->] (4,0) to (5,0); \node at (4.5,0.5) {$T$}; 
       \node at (5.5,0) {...};
       \draw  [thick, ->] (6,0) to (7,0); \node at (6.5,0.5) {$T$}; 
       \draw(8,0) circle (1cm); \node at (8,0){$T^{k}A$};
       \draw  [thick, ->] (8.5,1) to (10,2); \node at (9.2,1.7) {$T$}; 
       \draw(11,2.5) circle (1cm); \node at (11,2.5){$T^{k+1}A$};
       \draw  [thick, ->] (11,1.5) to (11,0.5); \node at (11.5,1) {$T$}; 
       \node at (11,0.2) {.};  \node at (11,0.0) {.}; \node at (11,-0.2) {.};
       \draw  [thick, ->] (11,-0.4) to (11,-1.5); \node at (11.5,-1) {$T$}; 
       \draw(11,-2.5) circle (1cm); \node at (11,-2.5){$T^{k+n}A$};
       \draw(8,-1.3) ellipse (1.35cm and 0.9cm); \node at (8,-1.3){$T^{k+n+1}A$}; 
       \node at (8,-0.7){$B$};   \node at (3,-2){$B:=T^{k}A\cap T^{k+n+1}A$};
       \draw  [thick, ->] (10,-2.4) to (8.9,-2); \node at (9.6,-1.6) {$T$}; 
\end{tikzpicture}
\end{center}
\caption{Intersection of images}\label{f:intersection}\end{figure}

As we already noted in Section~\ref{s:intro}, in the inhomogeneous
case (i.e. if $Q_s^t$ does depend on $s$) the direct statement in
Theorem~\ref{t:main} remains correct, but the inverse one fails.
Namely, it is possible that the sum in (\ref{e:main}) is infinite
for all $s$, but
$$\limsup_{t\to\infty}P(\xi_{s+t}\in A|\xi_s\in A)=0 ,$$
for some $s$, i.e. there is no recurrence. The reason for this is
that despite the ``chain of (pre)images'' of the set $A$
inevitably intersects itself an infinite number of times, the
original set needs not to be included to the intersections 
(see Fig.~\ref{f:intersection}). 

A counterexample may be constructed as a finite state discrete time
Markov chain with 4 states, whose graph of transitions is given by
the following diagram: %
\beq{e:diag}{ x_3\gets x_3\gets x_0\gets x_1\longleftrightarrow x_2 .}%
To be more precise, we consider a deterministic version of the
process.

\begin{example} Let $X:=\{x_i\}_{i=0}^3$, $x_i\ne x_j$ for $i\ne j$,
and we have a family of one-to-one maps $T_s:X\to X$ indexed by
$s\in\IZ$. Denote the only image of a state $x\in X$ at time $s$
by $T_s^1x$ and its unique pre-image by $T_s^{-1}x$. Then we set
$$ A:=x_0 = T_1^{-1}x_1,~~ T_2^{-1}x_1=x_2,~~T_3^{-1}x_2=x_1,~~
   T_1^1x_0=x_3, ~~ T_i^1x_3=x_3~\forall i .$$
\end{example}

Clearly the graph of transitions defined by this construction
satisfies the diagram (\ref{e:diag}). This shows that for the
measure $m$ uniformly distributed on $X$ and the state $x_0$ the
relation (\ref{e:main+}) holds true, but this state is
non-recurrent.

\subsection{Conservativity}

Let us discuss the connection of a somewhat unusual property
(\ref{e:main}) with more ``classical'' notions related to the
conservativity of a system.

First of all using the forward image of a set consider a version
of (\ref{e:main}) in the opposite time direction: %
\beq{e:main+}{\sum_{n\ge0} m(Q^{n}(A) \cap A) = \infty~~~
              \forall A\in\cB:~m(A)>0} %
and a very similar property formulated in terms of the action on
measures rather than on sets: %
\beq{e:main-mes}{\sum_{n\ge1} Q^{n}m(A) = \infty ~~~\forall
A\in\cB:~m(A)>0. } %
Consider also two more properties often associated with the
(forward/backward) {\em conservativity} in dynamical systems
theory: %
\beq{e:cons}{Q^1(A)\subseteq A\in\cB \Longrightarrow
                   m(A\setminus Q^1(A))=0,} %
\beq{e:cons-}{Q^{-1}(A)\subseteq A\in\cB \Longrightarrow
                   m(A\setminus Q^{-1}(A))=0.} %
In the first of these properties we either assume that
$Q^1(\cB)\subseteq\cB$ (i.e. 1-images of all measurable sets are
measurable), or have to consider the inner measure $m_*$ instead
of $m$ in the right hand side of the implication.

\begin{theorem}\label{t:X4}
(a) (\ref{e:main}) is equivalent to (\ref{e:cons}), ~~
(b) (\ref{e:main+}) does not imply PTP($m$), \\
(b') PTP($m$) does not imply (\ref{e:main-mes}), ~~
(c) (\ref{e:cons-}) implies (\ref{e:cons}). 
\end{theorem}

\proof We start with the direct statement in (a). Let
$Q^1(A)\subseteq A\in\cB$ and $m(A)>0$ (otherwise the statement is
obvious). Assume that $m(A) > m(Q^1(A))$. Then $m(B:=A\setminus
Q^1(A))>0$. Therefore
$$ \sum_{n\ge0}m(Q^{-n}B)=\infty $$
by (\ref{e:main}). On the other hand, by the construction
$$ Q^k(B) \cap B = \emptyset~~\forall k\in\IZ_+ .$$
Hence $Q^{-k}(B)\cap Q^{-n}(B) = \emptyset$ if $k\ne n$. Thus
$$ m(X) \ge m(\cap_{n\ge0}Q^{-n}(B)) = \sum_{n\ge0}m(Q^{-n}(B)) = \infty $$
by (\ref{e:main}). We came to the contradiction. %

It remains to consider the case $Q^1(A)\setminus A\ne\emptyset$.
In this case we set $B:=A$ and repeat the previous argument.

To prove (a) in the inverse direction consider the case $m(A)>0$,
$Q^1(A)\subseteq A$ and $m(A\setminus Q^1(A))=0$. Observe that
$Q^1(A)\subseteq A$ implies that $Q^{-1}(A) \supseteq A$ and hence
$Q^{-n}(A) \supseteq A~~\forall n\in\IZ_+$. Therefore
$$ \sum_{n\ge0}m(Q^{-n}(A)) \ge \sum_{n\ge0}m(A) = \infty ,$$
which proves validity of (\ref{e:main}). %

To demonstrate (b) consider the following example.

\begin{example}\label{ex1}
Let a Markov chain be defined on a compact space $X$ equipped with
a finite reference measure $m$ and divided into two disjoint parts
$X_1$ and $X_2$ such that $m(X_i)>0~\forall i$. The one-step
transition probabilities are defined by the relation:
$$ Q(x,B) := m(B\cap X_2)/m(X_2) .$$
\end{example}

Then for each set $A\subseteq X_1$ with $m(A)>0$ the property
(\ref{e:main+}) holds true, while no point from the set $A$ can
return back by dynamics and hence PTP($m$) is violated.

The claim (b') will be proven in Section~\ref{s:gen} during the
analysis of Example~\ref{ex2} (see Lemma~\ref{l:l2}).

(c) The property (\ref{e:cons-}) implies that for each backward
invariant set the restriction of the $m$-measure to the set $A$ is
preserved under dynamics. Therefore for such sets PTP($m$) is
obviously satisfied. On the other hand, from the proof of
Theorem~\ref{t:main} it follows that the assumptions
(\ref{e:main}) and $m(A)>0$ taken together imply the existence 
of a backward invariant set of positive $m$-measure. \qed 

Observe that neither of the properties $Q^1(A)\subseteq A$ and
$Q^{-1}(A)\subseteq A$ imply another one.

\section{Generalizations} \label{s:gen}

\subsection{Strong recurrence} As we already noted in the
Introduction, comparing our main results with their deterministic
counterparts or with results on recurrence and transience well
studied for the case of countable Markov chains (see e.g.
\cite{Fe,MT,Or,Sp}) one is tempted to make considerably stronger
statements assuming that the corresponding events take place with
probability one (instead of just being positive). To be precise
let us give a definition.

\bdef{ A point $x\in X$ is called {\em strongly Poincare
recurrent} with respect to a $\cB$-measurable set $A\ni x$ if
there exists an (arbitrary large) $t=t(x,A)$ such that
$Q^t(x,A)=1$. }

Consider an example.

\begin{example}\label{ex2}
Let a Markov chain be defined on a compact space $X$ equipped with
a finite reference measure $m$ and divided into two disjoint parts
$X_1$ and $X_2$ such that $m(X_i)>0~\forall i$. The one-step
transition probabilities are defined by the relation:
$$ Q(x,B) := \function{\frac{m(B\cap X_1)}{2m(X_1)}
                      + \frac{m(B\cap X_2)}{2m(X_2)}  &\mbox{if } x \in X_1 \\
                       \frac{m(B\cap X_2)}{m(X_2)} &\mbox{if } x \in X_2.}$$
\end{example}

The idea under this example is that the action of $Q$ on measures
splits them into two parts: one part remains confined in $X_1$,
while another part is moved to $X_2$ and never returns back.
Observe that without some special assumptions of mixing type this
situation is more or less generic.

\begin{lemma} Let a measurable set $A$ belongs to $X_i$ for some
$i\in\{1,2\}$, then each point $x\in A$ is Poincare recurrent, but
for any set $B\subseteq X_1$ there are no strongly Poincare
recurrent points.
\end{lemma}
\proof The first part of the claim follows from the fact that if
$A\subseteq X_1$, then
$$ Q(x,A) = \frac{m(A)}{2m(X_1)} > 0 ,$$
while if $A\subseteq X_2$, then
$$ Q(x,A) = \frac{m(A)}{m(X_2)} > 0 .$$
The second part is a consequence of the fact  that no point from the set
$X_1$ can return back under dynamics with probability one. \qed

\begin{corollary} In the example~\ref{ex2} PRP($m$) holds true,
but the similar statement for the strong recurrence fails.
\end{corollary}

\begin{lemma}\label{l:l2}
For a given measurable set $A\subset X_1$ with $m(A)>0$
$$ \sum_{n\ge1} Q^nm(A) < \infty .$$
\end{lemma}

The proof follows from the direct computation. This result
finishes the proof of part (b') of Theorem~\ref{t:X4} since it
demonstrates that PRP($m$) does not imply (\ref{e:main-mes}).

\bigskip

Stronger versions of other types of recurrence may be formulated
similarly to the Poincare recurrence.

Theorem~\ref{t:recur-t} claims that under very mild assumptions
the set of (weakly) recurrent points is of full reference measure.
On the other hand, for the Markov chain in Example~\ref{ex2} the
$m$-measure of the set of strongly recurrent points (both
topological and metric) is equal to $m(X_2)<m(X)$.

\subsection{Recurrence in the absence of invariant measures}

Let us answer the question, what happens to recurrence when the
corresponding Markov chain has no invariant measures. The simplest
example of this type has a deterministic nature and we assume the 
standard Borel $\sigma$-algebra and the corresponding topology.

\begin{example}\label{ex5} Let $X:=[0,1]$ and %
\beq{e:ex5}{
  Tx:=\function{x^2 &\mbox{if~}~~x >0 \\  1 &\mbox{otherwise}.} }
\end{example}

It is easy to see that the dynamical system defined by the
Example~\ref{ex5} has no invariant measures. We refer the reader
to \cite{Bl}, where more interesting examples of systems without
invariant measures are studied from the point of view of Birkhoff
type averaging.

\begin{lemma} Let $m$ be the Lebesgue measure on $X:=[0,1]$.
Then for each set $A$ containing a neighborhood of the origin all
points in $A$ are Poincare recurrent. The set of (both topological
and metrical) recurrent points coincides with the only point at
the origin.
\end{lemma}

In this example despite the absence of invariant measures 
the recurrent point is still present. To demonstrate that the 
situation may be even worse we use the map from Example~\ref{ex5} 
as a building block in the following construction.

\begin{example}\label{ex6} Let $X:=[0,1]$ and %
\beq{e:ex6}{
   Tx:=\function{4/5 &\mbox{if~}~~x=0 \\ 
                 x^2 &\mbox{if~}~~0<x<1/2 \\  
                 1-T(1-x) &\mbox{otherwise}.} }
\end{example}

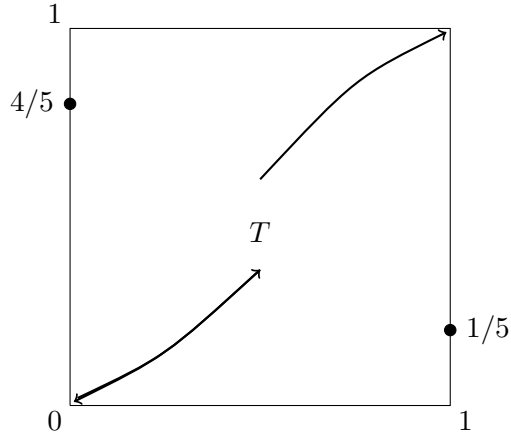
\begin{figure}
\begin{center}
\begin{tikzpicture} \small 
\put(0,0){\draw  (0, 0) to (5,0) to (5,5) to (0,5) to (0,0);
              \draw [thick, <-]  (0.05,0.05) .. controls (1.25,0.65) .. (2.5,1.8);  
              \draw [thick, ->]  (0.1,0.1) .. controls (1.25,0.65) .. (2.5,1.8);
              \draw [thick, ->]  (2.5,3.0) .. controls (3.75,4.35) .. (4.95,4.95);
              \draw [fill] (0,4.0) circle (0.075cm);  \draw (-0.5,4.0) node{$4/5$};
              \draw [fill] (5,1.0) circle (0.075cm);  \draw (5.5,1.0) node{$1/5$};
              \draw (-0.2,-0.2) node{$0$}; \draw (5.2,-0.2) node{$1$};
              \draw (-0.2,5.2) node{$1$};  \draw (2.5,2.3) node{$T$};
             }
\end{tikzpicture}
\end{center}
\caption{Map with no recurrent points}\label{f:norecurrence}\end{figure}

The graph of this map is depicted on Fig.~\ref{f:norecurrence}. 

\begin{lemma} The map $T$ has no invariant measures and 
no topological or metrical recurrent points.
\end{lemma}
\proof The 1st claim is an easy consequence of the fact that 
all points in $X$ converge under dynamics either to $0$ or to $1$. 
Thus the only candidates for the role of recurrent points are exactly 
these two points, however under dynamics each of them is mapped 
to the basin of attraction of the different point. \qed

\subsection{Continuous time}

Consider now a continuous time version of Theorem~\ref{t:main}.

\begin{theorem}\label{t:main-c}
Let $m$ be a finite measure on $(X,\cB)$ and let $Q_s^t$ does not
depend on $s$. Then PRP($m$) is equivalent to the existence of a
constant $\gamma>0$ such that %
\beq{e:main-c}{\sum_{n\ge1} m(Q^{-n\gamma}(A) \cap A) = \infty 
               ~~~\forall A\in\cB:~m(A)>0 .} %
\end{theorem}

The proof of this result follows from the time discretization with
time step $\gamma$ of the continuous Markov chain under study.
This gives a discrete time Markov chain $\tilde\xi_n$ with
one-step transition probabilities %
$$ \tilde Q(\cdot,\cdot):=Q^\gamma(\cdot,\cdot) .$$
By Theorem~\ref{t:main} the property PRP($m$) is equivalent to
$$ \sum_{n\ge1} m(\tilde Q^{-n}(A) \cap A) = \infty ~~~\forall
                    A\in\cB:~m(A)>0 ,$$ %
which immediately implies the claim of Theorem~\ref{t:main-c}.

\bigskip

In Theorem~\ref{t:main-c} we assume the lattice structure of the 
set of moments of time $\{n\gamma\}_{n\in\IZ}$ at which we observe 
our Markov chain. Can we use instead an arbitrary sequence of 
moments of time $t_1<t_2<\dots< t_n<\dots$ growing to infinity? 
Of course, one can make the time discretization with respect 
to these moments, but the corresponding diescrete time Markov 
chain will no longer be homogeneous.  

It is worth noting that since the proof of Theorem~\ref{t:recur-t}
does not use the discrete structure of time, it remains valid in
the continuous time setting as well.

\subsection{Multiple recurrence}

Furstenberg Theorem \cite{Fu} asserts that if $T$ is a measure
preserving bijection of a measure space $(X,\cB,\mu)$ with $\mu(X)
< \infty$ and if $A\in\cB$ with $\mu(A)>0$ and $k\ge2$ is any
integer, then there exists $n$ with
$$ \mu(A \cap T^nA \cap T^{2n}A \dots T^{(k-1)n}A ) > 0.$$
For $k = 2$, this is the Poincare recurrence property, and follows
readily from the fact that $T$ is measure preserving and
$\mu(X)<\infty$. For $k>2$ the proof of this result is quite nontrivial 
and despite of a large number of publications (see e.g. \cite{Kra} 
and further references therein) claiming that
this result holds true for a general measure preserving dynamical
system (rather than a bijection as in the original Furstenberg
publication), we were able to find a complete proof only for the
invertible case. No results in this direction are known for
general Markov processes.

Nevertheless we expect that a probabilistic version of this result
should hold true under the assumptions of Theorem~\ref{t:main},
namely that for any measurable set $A\in\cB$ with $m(A)>0$ and any
integer $k\ge2$ there exists a positive integer $n$ such that
$$ m(\cap_{j=0}^{k-1}Q^{-jn}(A)) > 0 .$$
Since the methods used to prove the original result due to
Furstenberg (which were further developed by his numerous
followers) are not applicable in the context of general Markov
chains (in particular, due to the non-invertibility of the
dynamics) a different approach needs to be developed. It is worth
noting that the original proof is quite involved and complicated.



\begin{thebibliography}{99}

\bibitem{Bl} Blank~M.
   Ergodic averaging with and without invariant measures. %
   Nonlinearity 30:12(2017), 4649-4664. arxiv:1709.06327 [math.DS]

\bibitem{Si} Cornfeld~I.P., Fomin~S.V., Sinai~Y.G.
  Ergodic Theory. New York: Springer-Verlag. 1982. 

\bibitem{Fe} Feller~W.
  An introduction to probability theory and its applications.
  V.1, Wiley. 1966

\bibitem{Fu} Furstenberg~H.
  Ergodic behavior of diagonal measures and a theorem of Szemeredi
  on arithmetic progressions. 
  Anal. Math. (1977) 31:1, 204-256.

\bibitem{Ha} Harris~T.E.
 The existence of stationary measures for certain Markov processes. 
  Matematika, 1960, Volume 4, Issue 1, 131-143.

\bibitem{Kra} Kra~B. 
  The Green-Tao Theorem on arithmetic progressions in the primes: an ergodic point of view. 
  Bull. Amer. Math. Soc. 43 (2006), 3-23.

\bibitem{MT} Meyn~S.P., Tweedie~R.L.
  Markov Chains and Stochastic Stability. Springer, New York, 1994.

\bibitem{Ni} Nitecki~Z. Differentiable dynamics. An introduction to the orbit
  structure of diffeomorphisms. Cambridge: The MIT Press. 1974

\bibitem{Or} Orey~S.
  Recurrent Markov chains.
  Pacific Journal of Mathematics, 9:3(1959), 806-827.

\bibitem{Sp} Spitzer~V.
  Principles of random walk. Nostrand. 1964

\end{thebibliography}
\end{document}